\documentclass[12pt,reqno]{amsart}
\usepackage[dvips]{graphicx}

\newtheorem{theorem}{{\bf Theorem}}[section]

\newtheorem*{proposition*}{{\bf Proposition}}

\newtheorem{lemma}[theorem]{{\bf Lemma}}
\newtheorem{lemma*}{{\bf Lemma}}

\newcommand{\CrossRatio}{{\bf c}}

\newcommand{\Uniform}{\mathit{u}}

\newcommand{\SLTwoC}{\mathrm{SL}_2(\mathbb{C})}
\newcommand{\PSLTwoC}{\mathrm{PSL}_2(\mathbb{C})}

\begin{document}

%
%

\title[ Circle Packings and Uniformization ]{%
  Circle Packings on Surfaces with Projective Structures
    and Uniformization }

\author[S. Kojima]{%
    Sadayoshi Kojima}
\address{%
        Department of Mathematical and Computing Sciences \\
        Tokyo Institute of Technology \\
        Ohokayama, Meguro \\
        Tokyo 152-8552 Japan}
\email{%
        sadayosi@is.titech.ac.jp}

\author[S. Mizushima]{%
    Shigeru Mizushima}
\address{%
        Department of Mathematical and Computing Sciences \\
        Tokyo Institute of Technology \\
        Ohokayama, Meguro \\
        Tokyo 152-8552 Japan}
\email{%
        mizusima@is.titech.ac.jp}

\author[S. P. Tan ]{%
    Ser Peow Tan}
\address{%
        Department of Mathematics \\
        National University of Singapore \\
        Singapore 117543, \\
        Singapore}
\email{%
        mattansp@nus.edu.sg}


\thanks{The third author gratefully acknowledge support
from the National University of Singapore academic research
grant R-146-000-031-112 and the hospitality of the Department of
Mathematical and Computing Sciences, Tokyo Institute of
Technology. }

%
%


\begin{abstract}
Let $\Sigma_g$  be a closed orientable surface of genus $g \geq 2$  and
$\tau$ a graph on $\Sigma_g$ with one vertex which lifts to a
triangulation of the universal cover.
We have shown that the cross ratio parameter space  $\mathcal{C}_{\tau}$
associated with  $\tau$,
which can be identified with the set of all pairs of
a projective structure and a circle packing on it
with nerve isotopic to  $\tau$,
is homeomorphic to ${\mathbb R}^{6g-6}$,
and moreover that
the forgetting map of  $\mathcal{C}_{\tau}$  to the space of
projective structures is injective.
In this paper, we  show that the composition of the forgetting
map with the uniformization
from  $\mathcal{C}_{\tau}$  to
the Teichm\"uller space ${\mathcal T}_g$ is proper.

\end{abstract}

\maketitle

%
%

\section{Introduction}
Let $\Sigma_g$ be a closed orientable surface of genus $g \ge 2$
and $\tau$ a graph on $\Sigma_g$ which lifts to a triangulation of
the universal cover $\widetilde{\Sigma}_g$.
In \cite{KMT},
we initiated the study of circle packings on $\Sigma_g$ with projective
structures with combinatorics controlled by $\tau$.
Since we discuss projective structures
which are complex structures in particular,
we orient  $\Sigma_g$  throughout this paper
so that the orientation is the same as the one coming from
the complex structure.

As the circle/disk is a fundamental
geometric notion for projective structures,
it makes sense to ask which surfaces admit a circle
packing with fixed nerve isotopic to  $\tau$,
and furthermore, if the circle
packings carried on these surfaces are rigid.
The main results of  \cite{KMT}  answer these two questions
when  $\tau$ has exactly one vertex.
In this case,
the cross ratio parameter space  $\mathcal{C}_{\tau}$  associated with
a graph  $\tau$  on  $\Sigma_g$,
which can be identified with the set of all pairs of
a projective structure and a circle packing on it
with nerve isotopic to  $\tau$,
is homeomorphic to  $\mathbb{R}^{6g-6}$.
Furthermore,
the forgetting map,
\begin{equation*}
    f : \mathcal{C}_{\tau} \to \mathcal{P}_g,
\end{equation*}
of  $\mathcal{C}_{\tau}$   to the space  $\mathcal{P}_g$  of
projective structures on  $\Sigma_g$  which forgets the packing is
injective. Namely, the packings are  in fact rigid.

On the other hand,
any projective structure on  $\Sigma_g$  has
a canonical underlying complex structure.
Thus assigning the underlying complex structure
to each projective structure,
we obtain the uniformization map
\begin{equation*}
    \Uniform : \mathcal{P}_g \to \mathcal{T}_g,
\end{equation*}
of  $\mathcal{P}_g$
to the Teichm\"uller space  $\mathcal{T}_g$
(thought of here as the space of complex structures).
By taking the Schwarzian derivative of the developing map,
a projective
structure can be identified with a holomorphic quadratic
differential over the underlying Riemann surface,
so the uniformization map is a complex vector bundle of
rank $3g-3$  over  $\mathcal{T}_g$.

For the genus one case, when $\tau$ has one vertex, it was shown
by Mizushima in \cite{Miz}, with slightly different language, that
the composition of the forgetting map with the uniformization map
is a homeomorphism. In \cite{KMT}, we conjectured that this holds
in general, regardless of the genus  $g \, (\geq 1)$ and the graph
$\tau$. The purpose of this paper is to take the first step
towards this conjecture by proving the following properness
theorem for  $\tau$ with one vertex:

\begin{theorem}\label{theo:maintheorem}
Let $\tau$ be a graph on $\Sigma_g$ ($g \ge 2$)
with one vertex which lifts to a triangulation of
$\widetilde{\Sigma}_g$ and $\mathcal{C}_{\tau}$
the cross ratio parameter space associated with  $\tau$.
Then the composition
$\Uniform \circ f : \mathcal{C}_{\tau} \to \mathcal{T}_g$  of
the forgetting map with the uniformization map is proper.
\end{theorem}

To complete the proof of the conjecture for such graphs $\tau$,
since  $\mathcal{C}_{\tau}$  in this case was shown to be
homeomorphic to  $\mathbb{R}^{6g-6}$,
it suffices to show that the map is locally injective.
This sort of question for the grafting map
based on Tanigawa's properness
theorem in  \cite{Tan} was settled by Scannell and Wolf in
\cite{SW}.
See also Faltings \cite{Fa} and McMullen \cite{McM} for
earlier proofs of special cases.
However,
it is not clear if the
proofs in the above cited papers can be extended to our setting.

The rest of this paper is organized as follows. In \S 2, we set
the notations,
and results required to prove the main theorem.
In \S 3,  we recall
the definition of the cross ratio parameter space
from \cite{KMT} and show that it is
properly embedded in the euclidean space.
In \S 4,
following the exposition in  \cite{KT},
we briefly review the Thurston coordinates of
$\mathcal{P}_g$  in terms of hyperbolic structures and
measured laminations,
and show that the projected image of  $f(\mathcal{C}_{\tau})$  to
the space of measured laminations on
$\Sigma_g$  is bounded.
The results up to that section are valid for any graph  $\tau$.
In \S 5,
under the assumption that  $\tau$  has exactly one vertex,
we show that the holonomy map from  $f(\mathcal{C}_{\tau})$
to the algebro-geometric quotient of the space of representations
of  $\pi_1(\Sigma_g)$  in  $\PSLTwoC$  up to conjugacy
is proper,
and deduce from this that the
projection of  $f(\mathcal{C}_{\tau})$  to
the space of the hyperbolic structures
in the Thurston coordinates is also proper.
Finally, in \S 6,
we complete the proof of the theorem using Tanigawa's
inequality in  \cite{Tan}.

%
%

\section{Preliminaries}

We first set the notation and recall results from \cite{KT, KMT}  to
be used in the rest of the paper.
Since $g$ and $\tau$ are fixed throughout this paper,
we will drop the suffices where possible to simplify the notation.
Let $\Sigma_g = \Sigma$ be a closed oriented surface of genus $g \ge 2$.

Attached to $\Sigma$,
we have the following spaces.
\begin{itemize}
\item The Teichm\"uller space ${\mathcal T}_g={\mathcal T}$
which parameterizes either the hyperbolic or the complex structures on
$\Sigma$, depending on the point of view.

\item The measured lamination space ${\mathcal M \mathcal
L}_g={\mathcal M \mathcal L}$ which is the completion of the space
of weighted simple closed curves on $\Sigma$.
It is known to
be homeomorphic to ${\mathbb R}^{6g-6}$.

\item The space ${\mathcal P}_g ={\mathcal P}$  of
projective structures on $\Sigma$,
namely a geometric structure locally modeled
on the Riemann sphere with transition functions in $\PSLTwoC$.

\item The algebro-geometric quotient  ${\mathcal X}_g={\mathcal X}$
of the space of representations of $\pi_1(\Sigma)$ in $\PSLTwoC$
up to conjugacy.
Taking the holonomy representation gives a map,
\begin{equation*}
    hol: {\mathcal P} \mapsto {\mathcal X},
\end{equation*}
which is a local homeomorphism, see \cite{Hej}.
\end{itemize}

Furthermore, we have the following spaces associated to $\Sigma$
and a graph  $\tau$ triangulating $\Sigma$
by the study in
\cite{KMT}.

\begin{itemize}
\item The cross ratio parameter space
$\mathcal{C}_{\tau}=\mathcal{C}$, which turns out to be a
semi-algebraic set in  ${\mathbb R}^{E}$, where  $E$  is the set
of edges of  $\tau$. An element  ${\bf c} \in \mathcal{C} \subset
\mathbb{R}^{E}$ is called a cross ratio parameter and will be
reviewed in the next section. It was shown in  \cite{KMT}  that
each cross ratio parameter determines a surface $S \in
\mathcal{P}$ with a projective structure together with a circle
packing  $P$ on $S$ with nerve isotopic to  $\tau$. Conversely,
each pair $(S, P)$ gives a cross ratio parameter  ${\bf c}$ in
$\mathcal{C}$.

\item The image  $\mathcal{I} = f(\mathcal{C})$  of  $\mathcal{C}$
by the forgetting map  $f : \mathcal{C} \to \mathcal{P}$, which is
the set of all projective structures on  $\Sigma$  which admit a
circle packing with nerve isotopic to  $\tau$. When  $\tau$  has
only one vertex, it was shown in  \cite{KMT}  that $f$  is a
homeomorphism onto the image $\mathcal{I}$  and hence  $f$
identifies  $\mathcal{C}$  with  $\mathcal{I}$ under this
condition.
\end{itemize}

Here is another view of  $\mathcal{P}$.

\begin{itemize}
\item
The Thurston coordinates of $\mathcal{P}$,
\begin{equation*}
    \mathcal{P} \cong \mathcal{T} \times \mathcal{ML},
\end{equation*}
which use the
pleated hyperbolic surface corresponding to a projective structure.
The first factor parameterizes the hyperbolic surface and the
second factor parameterizes the measured lamination coming from
the bending.
The coordinates are reviewed briefly in \S 4.
\end{itemize}

\medskip

%
%

\section{The Cross Ratio Parameter Space $\mathcal{C}$}

\begin{figure}[ht]
  \begin{center}
    \begin{picture}(300,140)
        \put(0,0){\scalebox{0.7}{\includegraphics{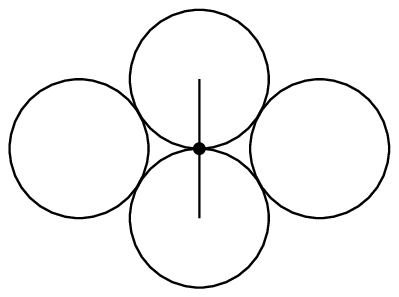}}}
        \put(160,0){\scalebox{0.7}{\includegraphics{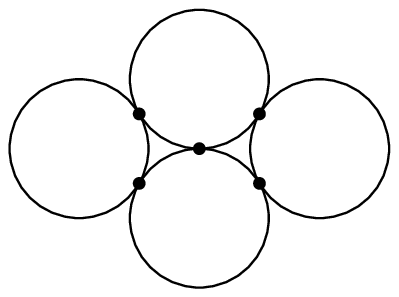}}}
\put(74,85){$\tilde{e}$}
\put(64,4){(a)} \put(191,81){$p_{23}$} \put(191,59){$p_{12}$}
\put(228,78){$p_{13}$} \put(257,81){$p_{34}$}
\put(257,58){$p_{14}$} \put(243,118){$C_3$} \put(293,50){$C_4$}
\put(246,20){$C_1$} \put(155,84){$C_2$} \put(224,4){(b)}
    \end{picture}
  \end{center}
  \caption{}
\label{Fig:Configuration}
\end{figure}

\begin{figure}[ht]
  \begin{center}
    \begin{picture}(140,160)
      \put(0,0){\scalebox{0.7}{\includegraphics{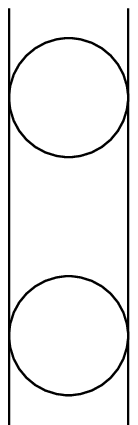}}}
\put(-28,109){$p_{14}=x\sqrt{-1}$}
\put(34,148){$C_1$}
\put(75,148){$C_3$} \put(85,109){$p_{34}$}
\put(127,143){$p_{13}=\infty$}
\put(85,28){$p_{23}=1$} \put(50,55){$C_2$}
\put(60,79){$C_4$}
\put(-3,28){$p_{12}=0$}
    \end{picture}
  \end{center}
  \caption[]{}
\label{Fig:CrossRatio}
\end{figure}

We first recall briefly how the cross ratio parameter space
$\mathcal{C}$ in  \cite{KMT} was defined. Suppose that $S$ is a
surface with projective structure which lies in  $\mathcal{I}
\subset \mathcal{P}$  and  $P$  is a circle packing on  $S$  with
nerve  $\tau$. To each edge $e$ of $\tau$ is associated a
configuration of four circles in the developed image surrounding a
preimage  $\tilde{e}$  of  $e$, see Figure
\ref{Fig:Configuration}(a). In \cite{KMT}, we  defined a cross
ratio of an edge  $e$ by taking the imaginary part of the cross
ratio of four contact points  $(p_{14},p_{23},p_{12},p_{13})$ of
the configuration chosen as in Figure \ref{Fig:Configuration}(b)
with orientation convention (for the definition of the cross ratio
of  four ordered points, see \cite{Ahl}). It is the modulus of the
rectangle obtained by normalizing the configuration by moving the
contact point $p_{13}$ to $\infty$, see Figure
\ref{Fig:CrossRatio}.

\begin{figure}[ht]
  \begin{center}
    \begin{picture}(140,120)
      \put(0,-20){\scalebox{0.7}{\includegraphics{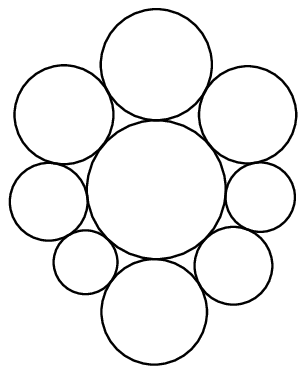}}}
\put(0,-20){ \put(66,67){$C$}
\put(66,110){$C_1$} \put(97,93){$C_2$} \put(99,65){$C_3$}
\put(91,41){$C_4$} \put(33,92){$C_m$} }
    \end{picture}
  \end{center}
  \caption[]{}
\label{Fig:Surround}
\end{figure}

Because the developing map is a
local homeomorphism, the cross ratios
of the edges must satisfy certain conditions, best expressed in
terms of the associated matrices.
If $e$ is an edge of $\tau$ with a cross ratio $x$,
we associate to  $e$  the matrix
$A={\left(\begin{array}{cc} 0 & 1 \\ -1 & x \\ \end{array}\right)}$.
Now,
if $v$ is a vertex of $\tau$ with valence $m$, we read off the edges
$e_1, \cdots, e_m$  incident to $v$  in a clockwise direction
to obtain a sequence of cross ratios  $x_1,\ldots,x_m$  associated to $v$.
Let
\begin{equation*}
    W_j=A_1 A_2 \cdots A_j
    ={\left(\begin{array}{cc} a_j & b_j \\ c_j & d_j \\
        \end{array}\right)}, ~~~j=1,\ldots, m
\end{equation*}
where $A_i$ is the matrix
${\left(\begin{array}{cc} 0 & 1 \\ -1 & x_i \\
\end{array}\right)}$ associated to $e_i$.
Then,
for each vertex $v$ of $\tau$, we have
\begin{eqnarray}
    W_v=A_1 A_2 \cdots A_m
    ={\left(\begin{array}{cc} -1 & 0 \\ 0 & -1 \\ \end{array}\right)} ,
\label{Eq:HolonomyCondition}
\end{eqnarray}
and
\begin{eqnarray}
    \left\{\begin{array}{l}
    a_j, \, c_j<0, \, b_j, \, d_j>0 \textrm{ for } 1\le j\le m-1\\
    \textrm{except for } a_1=d_{m-1}=0. \\
    \end{array}\right.
\label{Eq:SurroundOnceCondition}
\end{eqnarray}
Roughly, the first condition ensures the consistency of the
chain of circles surrounding the circle corresponding to $v$,
see Figure \ref{Fig:Surround},
and the second condition eliminates the
case where the chain surrounds the central circle more than once.
We note that it does not matter which edge
we start out with in the above.

Conversely, if an assignment to the edges of $\tau$,
\begin{equation*}
    {\bf c}: E \longrightarrow {\mathbb R},
\end{equation*}
satisfies conditions (\ref{Eq:HolonomyCondition})
and (\ref{Eq:SurroundOnceCondition}) for each vertex,
the map is a cross ratio parameter of some packing $P$
on a surface $S \in \mathcal{I}$, where
both $S$ and $P$ are determined by ${\bf c}$ (main lemma of \cite{KMT}).
Hence, the set
\begin{equation*}
    \mathcal{C} =
    \{ \CrossRatio : E \to {\mathbb R} \, \vert \,
    \CrossRatio \; \; \text{satisfies (1) and (2) for each vertex} \},
\end{equation*}
called the {\it cross ratio parameter space},
can be identified with the space of pairs $(S,P)$
of a surface $S$ with projective structure
and a circle packing $P$ on $S$,
and $\CrossRatio$ parameterizes the space of such pairs.

By the definition,
we see that $\mathcal{C}$ is a semi-algebraic set in
${\mathbb R}^{E}$ defined
by equations coming from (\ref{Eq:HolonomyCondition})
and inequalities coming from (\ref{Eq:SurroundOnceCondition})
for each vertex $v$.
We first show that
the strict inequalities of (\ref{Eq:SurroundOnceCondition}) can
be replaced by non-strict
inequalities, namely, for each vertex $v$ and with
the same notation as before,
consider the set of conditions
\begin{eqnarray}
    \left\{\begin{array}{l}
    a_j, \, c_j\le 0, \, b_j, \, d_j\ge 0 \textrm{ for } 1\le j\le m-1\\
    \textrm{except for } a_1=d_{m-1}=0 \\
    \end{array}\right.
\label{Eq:ExtendedSurroundOnceCondition}
\end{eqnarray}

\begin{lemma}\label{Lem:equivalentconditions}
The conditions
{\rm (\ref{Eq:HolonomyCondition}),(\ref{Eq:SurroundOnceCondition})}
$\Longleftrightarrow$
{\rm (\ref{Eq:HolonomyCondition}),(\ref{Eq:ExtendedSurroundOnceCondition})}.
\end{lemma}

\begin{proof}
It is sufficient to prove (\ref{Eq:HolonomyCondition}),
(\ref{Eq:ExtendedSurroundOnceCondition}) $\Longrightarrow$
(\ref{Eq:SurroundOnceCondition}).
We have the identity
\begin{eqnarray}
    \left(\begin{array}{cc} a_{j+1} & b_{j+1} \\
    c_{j+1} & d_{j+1} \end{array}\right) &=&
    \left(\begin{array}{cc} a_j & b_j \\ c_j & d_j
    \end{array}\right)
    \left(\begin{array}{cc} 0 & 1 \\ -1 & x_{j+1}
    \end{array}\right) \nonumber \\
    &=&
    \left(\begin{array}{cc} -b_j & a_j+b_j x_{j+1} \\
    -d_j & c_j+d_j x_{j+1}  \end{array}\right).
\label{Eq:productidentity}
\end{eqnarray}
Now  $a_2=-b_1=-1<0$, so by induction,
using (\ref{Eq:ExtendedSurroundOnceCondition}) and
(\ref{Eq:productidentity}), $a_{j+1}=-b_j<0$
for $j=1,\ldots, m-1$
since
$$a_2<0 \Longrightarrow b_2>0\Longrightarrow a_3<0
\Longrightarrow  \cdots \Longrightarrow b_{m-1}>0
\Longrightarrow a_m<0.$$
Similarly, since
$c_1=-1<0$, by induction, using (\ref{Eq:ExtendedSurroundOnceCondition}) and
(\ref{Eq:productidentity}),
$c_{j+1}=-d_j<0$   for $j=1, \ldots, m-2$ since
$$c_1<0 \Longrightarrow d_1>0\Longrightarrow c_2<0
\Longrightarrow  \cdots \Longrightarrow
d_{m-2}>0
\Longrightarrow c_{m-1}<0.$$
\end{proof}
This lemma tells us that the condition (\ref{Eq:SurroundOnceCondition})
does not divide a connected component of the algebraic set determined by
(\ref{Eq:HolonomyCondition}).
It just chooses appropriate connected components.
With this, we can easily prove

\begin{lemma}\label{Lem:proper}
The inclusion map of  $\mathcal{C}$  into  $\mathbb{R}^E$  is
proper,
where  $E$ is the set of edges of $\tau$.
\end{lemma}

\noindent {\it Remark.}
The above result holds for a general
$\tau$ with no restriction on the number of vertices.
\medskip

\begin{proof}
Let $\{\CrossRatio_n\}$ be a sequence of points in the
intersection of a compact set in  $\mathbb{R}^{E}$ with
$\mathcal{C}$. Then, there exists a convergent subsequence. Let
$\CrossRatio_\infty$ denote the limit of the subsequence. Each
$\CrossRatio_n$ satisfies the conditions
(\ref{Eq:HolonomyCondition}) and (\ref{Eq:SurroundOnceCondition}).
Each entry of the product of the matrices in the conditions
(\ref{Eq:HolonomyCondition}) and (\ref{Eq:SurroundOnceCondition})
is a polynomial in terms of ${x_i}$'s and thus is continuous.
Hence, the limit $\CrossRatio_\infty$ satisfies the conditions
(\ref{Eq:HolonomyCondition}) and
(\ref{Eq:ExtendedSurroundOnceCondition}), and also lies in
$\mathcal{C}$ by  lemma \ref{Lem:equivalentconditions}.
\end{proof}

%
%

\section{Image in  $\mathcal{ML}$  is Bounded}

Suppose  $S \in \mathcal{P}$  is a surface with projective
structure. Recall from \cite{KT} that associated to the developed
image $dev(\widetilde{S})$ is a collection of maximal disks
$\{D_i\}$  on  $\widetilde{S}$. To each maximal disk, there is a
set of at least two ideal boundary points which lie in the closure
of $D_i$ but not in  $\widetilde{S}$. Take the convex hull $V_i
\subset D_i$ of the ideal boundary points with respect to the
hyperbolic metric on the interior of  $D_i$. Then for two distinct
maximal disks $D_i$ and $D_j$, $V_i \cap V_j =\emptyset$ and the
collection  $\{ V_i \}$  gives a stratification of $\widetilde{S}$
equivariant under the action of the fundamental group. Note that
each  $V_i$  is either an ideal polygon or a line, except in the
trivial case where $S$ is a hyperbolic structure, in which case
there is only one maximal disk. Also every point $z$ in
$\widetilde{S}$ is contained in a unique convex hull $V_i$. We
call the maximal disk  $D_i$ containing $V_i$ a supporting maximal
disk of  $z$.

The polygonal parts of the stratification  $\{ V_i \}$  support
the hyperbolic metric.
Collapsing parallel lines in  $\{ V_i \}$  to a line and
quotienting the resultant by the action of the fundamental group,
we obtain a map of  $S$  to a hyperbolic surface  $H$.
This collapsing construction gives a map
\begin{equation*}
    \pi : \mathcal{P} \to \mathcal{T}
\end{equation*}
where  $\mathcal{T}$  is regarded here as the space of
hyperbolic structures on  $\Sigma$.

The stratification  $\{ V_i \}$  also defines a geodesic
lamination  $\lambda$  on  $H$  by taking the complement of the
interior of polygonal parts. Moreover, using the convex hull of
the ideal points of the maximal disk not in the disk but in the
3-dimensional hyperbolic space  $\mathbb{H}^3$, we can assign a
transverse bending measure supported on  $\lambda$. Hence we also
obtain a map
\begin{equation*}
    \beta : \mathcal{P} \to \mathcal{ML}.
\end{equation*}
Thurston showed  that the pair  $(\pi, \beta)$  of these maps is a
homeomorphism between  $\mathcal{P}$  and $\mathcal{T} \times
\mathcal{ML}$.

\begin{lemma}\label{Lem:compactbending}
$\beta(\mathcal{I})$  is bounded in  $\mathcal{M}\mathcal{L}$.
\end{lemma}

\noindent {\it Remark.}
The above result holds for a general
$\tau$ with no restriction on the number of vertices.

\begin{proof}
Let  $S$  be a surface
in  $\mathcal{I} = f(\mathcal{C}) \subset \mathcal{P}$,
$P$  a circle packing on  $S$  with nerve  $\tau$,
and  $H=\pi(S), \, \lambda=\beta(S)$.
The measured lamination  $\lambda$  can be pulled back
to a lamination with transverse measure
by blowing up the atomic leaves on  $H$  to parallel leaves in  $S$
with stretched transverse measure in a canonical way.
Hence we regard  $S$  as a surface with
this measured lamination  $\mu$.
To see boundedness of  $\beta(\mathcal{I})$,
it is sufficient to show that
the measure along each edge of  $\tau$  is uniformly bounded
since  $\tau$  generates the fundamental group of  $\Sigma$
and  $\mu$  collapses to  $\lambda$.

We can choose a reference point for each circle in  $P$  to
represent the vertex $v$ of $\tau$ such that the supporting
maximal disk of the point representing $v$ contains the circle and
take their preimage in $\widetilde{P}$ to get equivariant
reference vertices. Let $C_1$ and  $C_2$  be contact circles in
$\widetilde{P}$, and $D_i$ the supporting maximal disk of the
reference point  $v_i$ of $C_i$.   $D_i$ contains  $C_i$ for both
$i = 1, 2$, so these form a $\pi$ roof over the pleated hyperbolic
surface in $\mathbb{H}^3$ locally corresponding to
$dev(\widetilde{S})$, see Figure \ref{Fig:piroof}. So the total
transverse measure along a path $\tilde{e}$ between  $v_1$  and
$v_2$  contained in  $D_1 \cup D_2$  is bounded above by  $\pi$.
\begin{figure}[ht]
  \begin{center}
    \begin{picture}(246,86)
      \put(0,5){\scalebox{0.8}{\includegraphics{piroof2dim-2.eps}}}
      \put(30,43){$C_1$}
      \put(86,43){$C_2$}
      \put(10,8){$D_1$}
      \put(108,8){$D_2$}
      \put(160,0){\scalebox{0.3}{\includegraphics{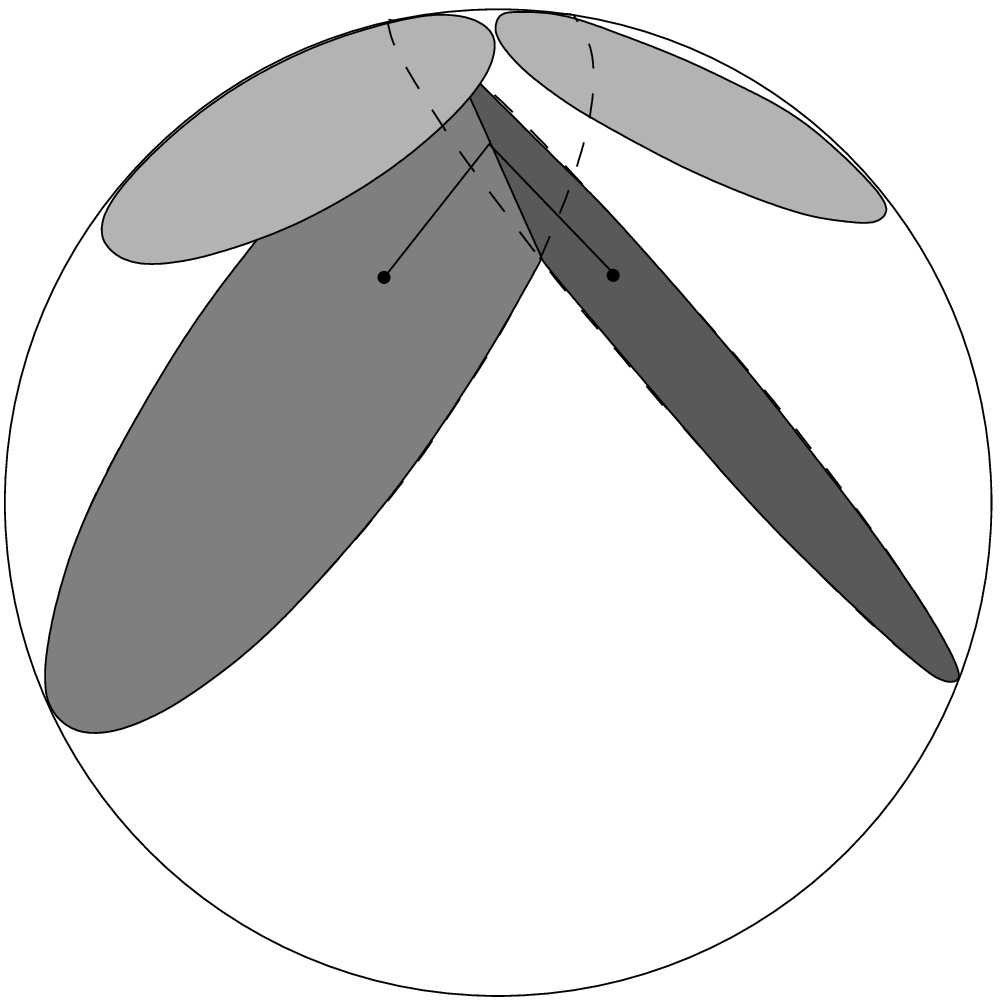}}}
    \end{picture}
  \end{center}
  \caption{}
  \label{Fig:piroof}
\end{figure}
%

%
\end{proof}

\medskip

%
%

\section{$hol$ and $\pi$ on Image are Proper}

In this section,
we restrict ourselves to the case
where $\tau$ has one vertex.
We first show

\begin{lemma}\label{Lem:holproper}
Suppose that  $\tau$  has only one vertex. Then there is a finite
subset $F$ of the fundamental group $\pi_1(\Sigma)$ such that if
$\{ {\bf c}_n \}$ is a sequence of points in  $\mathcal{C}$ which
escapes away from the compact sets, then there is an element $g
\in F$ and a subsequence of $\{ {\bf c}_n \}$ for which
$|tr(\rho_n(g))| \rightarrow \infty$, where  $\rho_n = hol(f({\bf
c}_n))$. In particular, the composition  $hol \circ f :
\mathcal{C} \to \mathcal{X}$  of the forgetting map with  $hol :
\mathcal{I} \to \mathcal{X}$  is proper.
\end{lemma}

\noindent {\it Remark.} Although the image of holonomy lies in
$\PSLTwoC$  and the trace makes sense only up to signs, its
absolute value is well-defined.

\begin{proof}
Since  $\tau$  has one vertex,
each edge $e_i$ of $\tau$ starts and ends
at the same vertex $v$,
so corresponds to a pair of elements
$g_i^{\pm 1} \in \pi_1(\Sigma)$,
and the set $\{g_i^{\pm 1}\}$
forms a generating set for $\pi_1(\Sigma)$.
$F$  is defined to be the set of
all words of length  $2$  in this generating set.

Passing to a subsequence, we may assume that the cross ratio ${\bf
c}_n(e)$  of some fixed edge $e$  of $\tau$ approaches $\infty$ as
$n \rightarrow \infty$. For each $S_n = f({\bf c}_n)$, consider
the developed image of $\widetilde{S}_n$ and in particular the
configuration of 6 circles $C_1,\ldots,C_6$ in $\widetilde{P}_n$
with corresponding vertices $v_1,\ldots,v_6$ of $\tilde \tau$ as
given in Figure \ref{Fig:TraceInfinityConfiguration}a. Here, the
edge $e_3=v_1v_3$ is the one whose cross ratio approaches $\infty$
as $n \rightarrow \infty$. Note that the configuration of  $C_j$
and  $v_j$  depends on  $n$.

To each  $n$, we may normalize the developed image so that it is
given by Figure \ref{Fig:TraceInfinityConfiguration}b, where the
tangency point $p_{13}$ between $C_1$ and $C_3$ is $\sqrt{-1}$.
The concatenation of the two directed edges $e_2=v_2v_1$ and
$e_4=v_1v_4$  which are the neighboring edges of $e_3$ about the
vertex $v_1$  corresponds to an element $g=g_4g_2^{-1}\in F$. Its
holonomy image $\varphi_n:= \rho_n(g) = \rho_n(g_4g_2^{-1})$ is an
element of $\PSLTwoC$ mapping $C_2$ to $C_4$. Note that in the
normalized picture, the radius of $C_4$ approaches zero as $n
\rightarrow \infty$ since ${\bf c}_n(e_3)\rightarrow \infty$.
Hence we can represent the images of  $0$, $1$ and $\infty$, all
lying on $C_2$, under $\varphi_n$ by $\sqrt{-1}+\varepsilon_1$,
$\sqrt{-1}+\varepsilon_2$ and $\sqrt{-1}+\varepsilon_3$
respectively, where $\varepsilon_1$, $\varepsilon_2$,
$\varepsilon_3\rightarrow 0$  as $n \rightarrow \infty$. Letting
$K_n=\frac{\varepsilon_2-\varepsilon_3}{\varepsilon_2-\varepsilon_1}$,
$\varphi_n^{-1}$ is then described by
\begin{equation*}
    \varphi_n^{-1}: z \mapsto
    K_n\frac{z-(\sqrt{-1}+\varepsilon_1)}{z-(\sqrt{-1}+\varepsilon_3)},
    ~~~ \hbox{and}
\end{equation*}
\begin{equation*}
    M_n=\frac{1}{\sqrt{K_n(\varepsilon_1-\varepsilon_3)}}
    \left(\begin{array}{cc}{K_n} & {-K_n(\sqrt{-1}
    +\varepsilon_1)} \\ {1} & {-(\sqrt{-1} +\varepsilon_3)} \\
    \end{array}\right)\in \SLTwoC
\end{equation*}
is a matrix representative of $\varphi_n^{-1}$.
We have
\begin{equation*}
    tr(M_n)=\frac{1}{\sqrt{\varepsilon_1-\varepsilon_3}}
    (\sqrt{K_n}-\frac{1}{\sqrt{K_n}}
    (\sqrt{-1}+\varepsilon_3))
\end{equation*}
and  $|tr(\varphi_n)|=|tr(\varphi_n^{-1})|=|tr(M_n)|$.
If $|tr(M_n)| \rightarrow \infty$,
we are done,
otherwise, we must have
$(\sqrt{K_n}-\frac{1}{\sqrt{K_n}}(\sqrt{-1}+\varepsilon_3))\rightarrow
0$ since $(\varepsilon_1-\varepsilon_3)\rightarrow 0$ as
$n \rightarrow \infty$.
This implies that $K_n$ approaches $\sqrt{-1}$
since $\varepsilon_3 \rightarrow 0$ too.
Geometrically, this means
that the angle formed by $\sqrt{-1}+\varepsilon_1$,
$\sqrt{-1}+\varepsilon_2$ and $\sqrt{-1}+\varepsilon_3$ approaches
$\pi/2$, and hence $\sqrt{-1}+\varepsilon_1$ and
$\sqrt{-1}+\varepsilon_3$ approach diametrically opposite
positions on the circle $C_4$.
In other words, the arc
$\alpha_1=\varphi_n ([-\infty,0])$ occupies  half of the circle $C_4$
in the limit.

On the other hand, the inferior arc $\alpha_2$ on
$C_4$ connecting $p_{14}$ and $p_{34}$ also occupies half of $C_4$
in the limit since
the radius of  $C_4$,  which we denote by  $rad(C_4)$,
approach  $0$,
where  $p_{jk}$ is the
tangency point between $C_j$ and $C_k$.
\begin{figure}[ht]
  \begin{center}
    \begin{picture}(280,110)
      \put(0,15){\scalebox{0.8}{\includegraphics{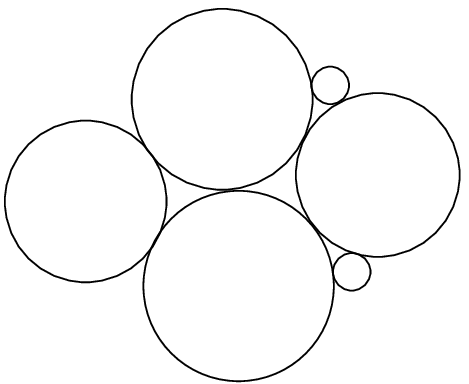}}}
      \put(50,0){(a)}
      \put(20,50){$p_{12}$}
      \put(19,65){$p_{23}$}
      \put(48,51){$p_{13}$}
      \put(74,53){$p_{14}$}
      \put(72,68){$p_{34}$}
      \put(-13,53){$C_2$}
      \put(65,98){$C_3$}
      \put(108,58){$C_4$}
      \put(65,10){$C_1$}
      \put(104,91){\vector(-4,-1){20}}
      \put(107,88){$C_5$}
      \put(108,34){\vector(-4,1){20}}
      \put(111,29){$C_6$}
      \put(140,20){\scalebox{0.8}{\includegraphics{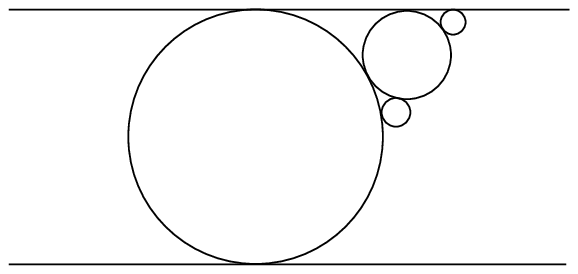}}}
      \put(197,0){(b)}
      \put(196,19){$p_{12}$}
      \put(214,17){$=0$}
      \put(196,92){$p_{13}$}
      \put(201,100){$\rotatebox{90}{=}$}
      \put(190,112){$\sqrt{-1}$}
      \put(210,70){$p_{14}$}
      \put(233,92){$p_{34}$}
      \put(230,40){$C_1$}
      \put(147,32){$C_2$}
      \put(147,90){$C_3$}
      \put(248,66){$C_4$}
      \put(273,78){\vector(-4,1){20}}
      \put(276,73){$C_5$}
      \put(259,50){\vector(-2,1){20}}
      \put(261,44){$C_6$}
    \end{picture}
  \end{center}
  \caption{}
  \label{Fig:TraceInfinityConfiguration}
\end{figure}

\noindent
{\bf Claim:}  $\alpha_1$ and $\alpha_2$ are the images under
$\varphi_n$ of two distinct, non-adjacent segments of $C_2$ and
thus are distinct non-adjacent segments on $C_4$.
\medskip

\noindent
{\it Proof of Claim.}
It suffices to show that the
image of $[-\infty,0]$ on $C_1$ by $\rho_n^{-1}(g_2)$ and the image
of the arc $\alpha_2$ by $\rho_n^{-1}(g_4)$ are non-adjacent
segments of $C_1$. The pair of edges $(e_2,e_3)$ are either part
of an intersecting or non-intersecting triple of edges, since they
are adjacent (see \cite{KMT}, \S 4).
Similarly, the pair of edges
$(e_3,e_4)$ are also part of an intersecting or non-intersecting
triple.
For all possibilities, $\rho_n^{-1}(g_2)[-\infty, 0]$ and
$\rho_n^{-1}(g_4) (\alpha_2)$ are non-adjacent segments of $C_1$,
see Figure \ref{Fig:nonadjacentsegments} (a),(b) and (c).
\qed
\medskip

\begin{figure}[ht]
  \begin{center}
    \begin{picture}(270,235)
      \put(0,140){
        \put(0,0){\scalebox{0.65}{\includegraphics{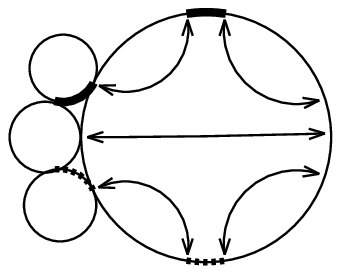}}}
        \put(3,76){$C_4$}
        \put(-5,48){$C_3$}
        \put(0,20){$C_2$}
        \put(110,22){$C_1$}
        \put(24,68){$\alpha_2$}
        \put(3,3){$[-\infty,0]$}
        \put(19,11){\vector(1,2){10}}
        \put(68,92){$\rho_n^{-1}(g_4)(\alpha_2)$}
        \put(65,-5){$\rho_n^{-1}(g_2)([-\infty,0])$}
        \put(0,-20){(a) Two non-intersecting}
        \put(18,-33){triples}
      }
      \put(150,140){
        \put(0,0){\scalebox{0.65}{\includegraphics{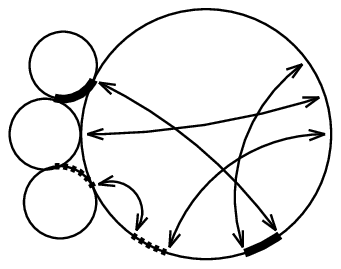}}}
        \put(3,76){$C_4$}
        \put(-5,48){$C_3$}
        \put(0,20){$C_2$}
        \put(74,90){$C_1$}
        \put(0,-20){(b) Intersecting and}
        \put(18,-33){non-intersecting triples}
      }
      \put(0,0){
        \put(0,0){\scalebox{0.65}{\includegraphics{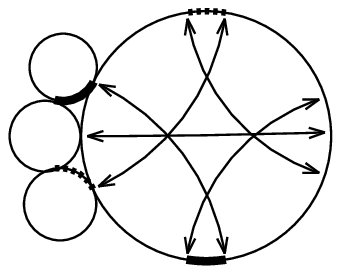}}}
        \put(3,76){$C_4$}
        \put(-5,48){$C_3$}
        \put(0,20){$C_2$}
        \put(115,39){$C_1$}
        \put(0,-10){(c) Two intersecting triples}
      }
    \end{picture}
  \end{center}
  \caption{}
  \label{Fig:nonadjacentsegments}
\end{figure}
%

From the claim, the subtended angle between $p_{45}$ and $p_{34}$
on $C_4$ in Figure \ref {Fig:TraceInfinityConfiguration} (b)
approaches zero, as does the angle between $p_{46}$ and $p_{14}$.
Thus, the ratios of the radius of $C_5$ and $C_6$ to that of
$C_4$, $\frac{rad(C_5)}{rad(C_4)}$ and $\frac{rad(C_6)}{rad(C_4)}$
approach zero and the cross ratios of the two edges $v_1v_4$ and
$v_3v_4$ both approach $\infty$. Repeating the argument with the
roles of $C_2$ and $C_4$ reversed, we see that either
$|tr(\varphi_n)|\rightarrow \infty$ or the cross ratios of
$v_1v_2$ and $v_2v_3$ both approach $\infty$. In other words,
either $|tr(\varphi_n)|\rightarrow \infty$ or the cross ratios of
all neighbor edges approach $\infty$. By induction, using the edge
$e_4$ instead of $e_3$ and repeating the above argument, we see
that either some element $g \in F$ has holonomy with diverging
trace, or the cross ratios of all the edges of $\tau$ diverge to
$\infty$. However, in the latter case, the $(2,2)$ term  $d_m$  of
the matrix corresponding to the word $W_m$ defined in \S 3  has a
dominating term $x_1 x_2 \cdots x_{m}$ and diverges to $\infty$.
This contradicts the fact that $d_{m}=-1$ for all points of the
sequence by condition (\ref{Eq:HolonomyCondition}) in \S3. This
completes the proof.
\end{proof}

Lemma \ref{Lem:holproper} implies

\begin{lemma}\label{Lem:p1proper}
Suppose that  $\tau$  has only one vertex.
The composition  $\pi \circ f : \mathcal{C} \to \mathcal{T}$
of the forgetting map with the collapsing map
$\pi : \mathcal{P} \to \mathcal{T}$  in
the Thurston construction is proper.
In particular,
the restriction of  $\pi$  to  $\mathcal{I}$  is proper.
\end{lemma}

\begin{proof}
Suppose that $\{ {\bf c}_n \}$ is a sequence in $\mathcal{C}$ which
escapes away from the compact sets.
By passing to a subsequence, we may assume by lemma
\ref{Lem:holproper} that there is an element $g \in \pi_1(\Sigma)$ with
holonomy $\varphi_n=\rho_n(g)$ for
which $|tr(\varphi_n)| \rightarrow \infty$.
$g$ corresponds to a closed
curve $\gamma$ on $\Sigma$ and the length $l_n(\gamma)$ of the
geodesic representative of $\gamma$ on the collapsed
surface $H_n=\pi(S_n)$ satisfies the inequality
\begin{equation*}
    l_n(\gamma) \ge d_{\mathbb{H}^3}(z_n, \varphi_n(z_n))
\end{equation*}
where $z_n$ is any point on the axis of $\varphi_n$  in  $\mathbb{H}^3$  and
$d_{\mathbb{H}^3}(z_n, \varphi_n(z_n))$ is the hyperbolic distance between
$z_n$ and $\varphi_n(z_n)$. Since $|tr(\varphi_n)|\rightarrow
\infty$, $d_{\mathbb{H}^3}(z_n, \varphi_n(z_n)) \rightarrow \infty$ and hence
$l_n(\gamma)\rightarrow\infty$.
It follows that $\{ H_n \}$  escapes away from the compact sets
in the Teichm\"uller space $\mathcal{T}$.
\end{proof}

%
%

\section{Proof of Theorem}

We recall the following result from \cite{Tan}, the notation
has been modified slightly to fit into our framework.

\begin{theorem}[\cite{Tan}, theorem 3.4]
\label{theo:tanigawa} Let  $S=(H, \lambda)$ be a surface
with projective structure
where $H=\pi(S)\in \mathcal{T}$ and $\lambda =\beta(S) \in
\mathcal{M}\mathcal{L}$  and let $X$ be the underlying
Riemann surface.
Let $h:X \rightarrow H$ denote the
harmonic map with respect to the hyperbolic metric on $H$ and
$\mathcal{E} (h)$ its energy. Then
\begin{equation}\label{eqn:tanigawa}
    l_H(\lambda) \le
    \frac{l_H(\lambda)^2}{E_X(\lambda)} \le 2\mathcal{E} (h) \le
    l_H(\lambda)+8\pi(g-1),
\end{equation}
where $l_H(\lambda)$ is the
hyperbolic length of $\lambda$ on $H$, and $E_X(\lambda)$ is the
extremal length of $\lambda$ on  $X$.
\end{theorem}

We are now ready to prove the theorem.

\begin{proof}[Proof of Theorem 1.1]
Assuming that  $\tau$  has one vertex,
we follow the argument of Tanigawa in \cite{Tan}.
Let  $\{ {\bf c}_n \}$  be a sequence of points in  $\mathcal{C}$
which escapes away from the compact sets,
$f({\bf c}_n) = S_n=(H_n, \lambda_n) \in \mathcal{T} \times \mathcal{ML}$
and  $X_n = \Uniform (S_n)$  the corresponding
underlying Riemann surface.
By lemma \ref{Lem:p1proper},
we may assume that $\{ H_n \}$  escapes away from the compact sets
in $\mathcal{T}$.
Furthermore, by lemma \ref{Lem:compactbending},
since $\{ \lambda_n \}$
lies in a compact subset of $ \mathcal{ML}$,
we may assume by taking a subsequence if necessary
that $\lambda_n \rightarrow \lambda$ for some fixed measured
lamination $\lambda$.  Let $h_n$ be the harmonic map from $X_n$ to
$H_n$. We consider the two cases
\begin{itemize}
\item[(i)] $\sup_n l_{H_n}(\lambda_n) < \infty$, or

\item[(ii)] $\lim_{n\rightarrow \infty}l_{H_n}(\lambda_n) =
\infty$.
\end{itemize}
First assume that case (i) holds.
If $\{X_n\}$ stays in a compact
subset of  $\mathcal{T}$, then since $\{ H_n \}$
escapes away from the compact sets,
by a result of M. Wolf \cite{Wol},
$\mathcal{E}(h_n)\rightarrow \infty$.
This contradicts (i) and the right inequality of
(\ref{eqn:tanigawa}).
Hence $\{ X_n \}$  escapes away from the
compact sets.

Next suppose that case (ii) holds. Then by (\ref{eqn:tanigawa})
\begin{equation*}
    \lim_{n\rightarrow\infty}
    E_{X_n}(\lambda_n)=\lim_{n\rightarrow\infty}
    E_{X_n}(\lambda)=\lim_{n\rightarrow\infty}
        (l_{H_n}(\lambda_n)+O(1))=\infty,
\end{equation*}
and so $\{X_n \}$  escapes away from the compact sets as well
and we are done.
\end{proof}

%
%

\end{document}